\documentclass[12pt]{article}
\usepackage{amsfonts, amsmath, amstext, amssymb}

\usepackage[dvips]{graphicx}
\usepackage{latexsym}
\usepackage{color}

\newtheorem{thm}{Theorem}[section]
\newtheorem{prop}[thm]{Proposition}

\newtheorem{lemma}[thm]{Lemma}

\newtheorem{cor}[thm]{Corollary}

\newtheorem{definitiontemp}[thm]{Definition}
\newenvironment{defn}{\begin{definitiontemp}
\normalfont}{\end{definitiontemp}}

\newcommand\aut[1]{\text{Aut}(#1)}

\def\bec{\begin{cor}}
\def\enc{\end{cor}}

\def\bet{\begin{thm}}
\def\ent{\end{thm}}

\def\becor{\begin{cor}}
\def\encor{\end{cor}}

\def\bel{\begin{lem}}
\def\enl{\end{lem}}

\def\bedef{\begin{defn}}
\def\endef{\end{defn}}

\def\bep{\begin{prop}}
\def\enp{\end{prop}}

\newenvironment{pf}{\begin{trivlist}\item[\hskip\labelsep
{\it Proof.}]}{\end{trivlist}}

\newcommand{\set}[2]{\ensuremath{ \{ #1 : #2 \} }}

\renewcommand{\deg}[1]{\ensuremath{\text{deg}(#1)}}

\newcommand{\Z}{\mathbb{Z}}
\newcommand{\F}{\mathbb{F}}
\newcommand{\Q}{\mathbb{Q}}
\newcommand{\R}{\mathbb{R}}

\newcommand{\A}{\mathfrak{A}}

\newcommand{\Qbar}{\overline{Q}}

\newcommand{\Ftilde}{\tilde{F}}

\newcommand{\ptilde}{\tilde{p}}
\newcommand{\qtilde}{\tilde{q}}

\newcommand{\ztilde}{\tilde{z}}

\newcommand{\B}{\mathfrak{B}}
\newcommand{\Bbar}{\overline{B}}

\newcommand{\mR}{\mathcal{R}}

\newcommand{\avec}{\vec{a}}
\newcommand{\bvec}{\vec{b}}

\newcommand{\xvec}{\vec{x}}

\newcommand{\zvec}{\vec{z}}

\newcommand{\la}{\langle}
\newcommand{\ra}{\rangle}
\def\embeds{\hookrightarrow}

\newcommand{\hid}{h_{\text{id}}}

\def\diverges{\!\uparrow}
\def\converges{\!\downarrow}
\newcommand{\at}{\char'100}

\newcommand{\qed}{\hbox to 0pt{}\nobreak\hfill\rule{2mm}{2mm}}

\newcommand{\dom}[1]{\text{dom}(#1)}

\newcommand{\bfc}{\boldsymbol{c}}
\newcommand{\bfd}{\boldsymbol{d}}
\def\bfz{\boldsymbol{0}}
\def\s01{\ensuremath{\Sigma^0_1}}
\def\d02{\ensuremath{\Delta^0_2}}
\def\phi{\varphi}

\def\res{\!\!\upharpoonright\!}

\begin{document}

\title{Computable Categoricity for Algebraic Fields with
Splitting Algorithms}
\author{
Russell Miller\thanks{Partially supported by Grant \# DMS -- 1001306 from
the National Science Foundation, by Grant \# 13397 from the
Templeton Foundation, by Queens College Research
Enhancement Program Award \# 90927-08 08,
and by grants numbered 61467-00 39,
62632-00 40, and 63286-00 41 from
The City University of New York PSC-CUNY Research Award Program.
} \&
Alexandra Shlapentokh\thanks{Partially supported by Grant \# DMS--0650927
from the National Science Foundation, by Grant \# 13419 from the
Templeton Foundation, and by an ECU Faculty Senate Summer 2011 Grant.}}
\maketitle

\begin{abstract}
A computably presented algebraic field $F$ has a \emph{splitting algorithm}
if it is decidable which polynomials in $F[X]$ are
irreducible there.  We prove that such a field is
computably categorical iff it is decidable which pairs
of elements of $F$ belong to the same orbit under
automorphisms.  We also show that this criterion is equivalent
to the relative computable categoricity of $F$.
%
\end{abstract}

\section{Introduction}
\label{sec:intro}

Computable categoricity is a basic concept in computable
model theory.  It holds of those computable structures $\A$
for which the classical concept of being isomorphic to $\A$
is equivalent to its analogue in the context of computable
structures and computable isomorphisms.
Thus it fits squarely into the program of computable
model theory, which studies how difficult it is
to give effective versions of classical
model-theoretic concepts and constructions.
(A summary of relevant definitions and standard results
appears in Sections \ref{sec:computability} and
\ref{sec:fields}.)

Assorted versions of this concept have been proposed.
The strongest, \emph{computable stability},
holds of a computable structure $\A$ iff every classical
isomorphism from $\A$ onto any other computable structure
is itself computable.  This holds, for instance,
for computable presentations of the structure
$(\Z, S)$, the integers under the successor function.
The most common version, \emph{computable categoricity},
which will be the focus of this paper, is not as
stringent:  a structure $\A$
is computably categorical iff, whenever there
exists a classical isomorphism $f$ from $\A$
onto a computable structure $\B$, there also exists
a computable isomorphism $g$ from $\A$ onto $\B$,
although $f$ itself may fail to be computable.
This version has been generalized to \emph{relative
computable categoricity}, which requires that
for every structure $\B$ which has domain $\omega$
and is isomorphic to the computable structure $\A$,
there must exist an isomorphism from $\A$ onto $\B$
which is computable in the Turing degree of the
(not necessarily computable) structure $\B$.
Finally, there are weaker versions such as
\emph{$\Delta^0_{n+1}$-categoricity}, and relativizations of these:
$\A$ is $\Delta^0_{n+1}$-categorical iff every
computable structure $\B$ isomorphic to $\A$
has an isomorphism onto $\A$ which is $\bfz^{(n)}$-computable.
Such versions essentially study how much information
about the structures is needed to compute isomorphisms.

Among these, computable categoricity remains
the most widely studied concept.  It is often equivalent
to relative computable categoricity, but exceptions are known
to exist; see \cite{K96} for an exception, and \cite{G75}
for conditions implying equivalence.  Traditionally, the main question
has been to determine, for a particular class of structures,
some structural criterion which is equivalent to computable
categoricity.  In early examples, from around 1980,
Dzgoev, Goncharov, and Remmel showed (independently;
see \cite{GD80, R81a})
that a linear order is computably categorical iff it
has only finitely many pairs of adjacent elements,
and Remmel also showed in \cite{R81b}
that a Boolean algebra is computably categorical
iff it has only finitely many atoms.  In both cases,
the structural criterion identifies the obstacle to computing
isomorphisms, in the class of structures under consideration.
On the other hand, the criteria equivalent to computable categoricity for trees
(viewed either as partial orders, or under the meet relation,
and also with distinguished subtrees),
established in \cite{LMMS05} and \cite{M05} by Lempp, McCoy,
Miller, and Solomon and in \cite{KKM05} by Kogabaev, Kudinov,
and Miller, are not easy to describe in any known way,
even though they are ``structural,'' in any reasonable sense of the word.
In terms of computational complexity, they are $\Sigma^0_3$,
as defined in Section \ref{sec:computability},
just like the conditions for linear orders and Boolean algebras.
Indeed, for all of these structures, computable categoricity
is equivalent to relative computable categoricity, and
it is known from work by Ash, Knight, Manasse, and Slaman
in \cite{AKMS89}, and independently by Chisholm in \cite{Chis},
that the computational complexity of relative computable categoricity
is always $\Sigma^0_3$.

Our intention in this paper is to give a criterion for computable
categoricity for algebraic fields with splitting algorithms. This should be viewed
as a first step towards a criterion for computable fields in general,
for which the question of computable categoricity has long
been studied and has proven highly intractable.  The basic
definitions regarding computable fields appear in
Section \ref{sec:fields}.  Apart from
algebraically closed fields, we believe that ours is the first
result to characterize computable categoricity for any natural
class of fields.  (An algebraically closed field is computably categorical
iff it has finite transcendence degree over its prime subfield.
The forwards direction was established by Ershov in
\cite{E77}; while the converse follows from the very first known
consideration of computable categoricity in the literature,
which forms part of the paper \cite{FS56} by Frohlich and Shepherdson.)
It is known that finite transcendence degree does not characterize
computable categoricity for fields in general:  in \cite{E77} Ershov
gave a field which is algebraic over its prime subfield,
yet is not computably categorical, and more recently Miller and Schoutens
disproved the converse in \cite{MS08}, by building
a computably categorical field of infinite transcendence degree over $\Q$.

A computable field $F$ is algebraic with a splitting algorithm if $F$
is an algebraic extension of its own prime subfield ($\Q$ or $\F_p$)
such that the set of reducible polynomials in $F[X]$ is decidable.
Our criterion for computable categoricity, which holds
on all such computable fields $F$, is that the orbit relation on the field
must be decidable.  In a certain sense, this criterion
is not quite as structural as previous ones for other
classes of structures; earlier criteria often used attributes,
such as finiteness, which are not expressible in first-order
model theory, but we are not aware of other classes
for which it is necessary to use the notion of decidability
to characterize computable categoricity.  Of course,
we do not prove here that it is actually necessary to do so
for algebraic fields with splitting algorithms either, but
the simplicity of the statement of our result, combined
with the lack of any other criteria after much study by
many researchers, leads us to believe empirically that
ours is as ``good'' a criterion as one is likely to find.
In complexity terms, the statement of our criterion is $\Sigma^0_3$,
just like those for many other classes, and it is equivalent to relative computable
categoricity, which would force it to be $\Sigma^0_3$
in any case.  Indeed, computable categoricity is readily
seen to be $\Sigma^0_3$-hard for such fields, and so
our criterion has the lowest possible computational complexity.

Recent work by Miller also studied categoricity for computable
algebraic fields, but using the weaker notion of $\bfd$-computable
categoricity.  A computable structure $\A$ is \emph{$\bfd$-computably categorical}
if every computable structure isomorphic to it is isomorphic via
some $\bfd$-computable isomorphism.  This says that the degree
$\bfd$ contains sufficiently much information to compute
an isomorphism whenever one exists.  Normally, for structures
which are not computably categorical, one finds
that they are $\bfz'$-computably categorical,
or $\bfz^{(m)}$-computably categorical for some $m$,
or possibly $\bfz^{(\alpha)}$-computably categorical
for some ordinal $\alpha$, and it is common that, for the least
such $m$, $\bfz^{(m)}$ gives a sharp lower bound for such degrees $\bfd$,
especially if we quantify over all structures in a class.
The surprising result in \cite{M09b} was that,
although not all computable algebraic fields are $\bfz'$-categorical,
there is a degree $\bfd$ with $\bfd'=\bfz''$ such that all
such fields are $\bfd$-computably categorical.
The unusual degree involved here suggests that
the question of computable categoricity for fields,
even just for algebraic fields, is somehow not as straightforward
as for many other structures.

On the other hand, the current work \cite{HKMS11}, by Hirschfeldt,
Kramer, and the present authors, shows that relative computable
categoricity does have a fairly reasonable structural characterization
for computable algebraic fields.  Their criterion does involve computability,
just as does ours in this paper, but it can be expressed in a generally
understandable way.  However, they also show that not all computably
categorical algebraic fields are relatively computably categorical,
and that the criterion from this paper, for algebraic fields
with splitting algorithms, fails to extend (at least in the obvious ways)
to computable algebraic fields without splitting algorithms.
Rather than producing a structural criterion, they show that the definition of
computable categoricity for algebraic fields is $\Pi^0_4$,
and that for such fields, the property of being computably
categorical is $\Pi^0_4$-complete.  Thus they establish that no structural criterion
can express computable categoricity in any simpler way
than the basic definition.  Of course, it is still possible that
some structural criterion of complexity $\Pi^0_4$ (or higher)
might characterize computable categoricity in a more transparent
way than the definition itself does; but in terms of computational complexity,
one cannot simplify the definition at all.

The next two sections of the paper contains definitions and background
on computability theory and on computable fields, along with a number of classical results
about fields which will be useful later on.  This much is sufficient for us
to introduce the problem, in Section \ref{sec:EX}, with
some basic cases of algebraic fields and computable categoricity.
Sections \ref{sec:isotrees} and \ref{sec:Galoisaction}
describe further concepts necessary for the main theorem:
the isomorphism tree, and the orbit relation.  Finally,
in Section \ref{sec:splitalg}, we prove the full result
for computable algebraic fields with splitting algorithms.

\section{Background in Computability}
\label{sec:computability}

We recall here the concepts from computability theory which will be
essential to our work on fields.  Computable
functions are defined in \cite{S87}, and indeed,
several very different definitions give rise to the
same class of functions.  Functions on the set $\omega$
of nonnegative integers are usually identified with
their graphs in $\omega^2$, and we then code $\omega^2$
into $\omega$, so that the graph corresponds to a
subset of $\omega$; conversely, for our purposes, a subset
of $\omega$ may be identified with its characteristic function.
The partial computable functions
(those for which the computation procedure halts
on certain inputs from $\omega$,
but not necessarily on all of them) can be enumerated
effectively, and are usually denoted as $\phi_0,\phi_1,\ldots$,
with the index $e$ coding the program for computing
$\phi_e(x)$ on $x\in\omega$.  The domains of these functions
constitute the \emph{computably enumerable sets}, and we write
$W_e$ for the domain of $\phi_e$.  These are precisely the sets
which are definable by $\Sigma^0_1$ formulas,
i.e.\ sets of the form
$$ \set{x\in\omega}{\exists y_1\cdots\exists y_m~R(x,y_1,\ldots,y_m)},$$
where $m\in\omega$ is arbitrary and $R$ may be
any computable subset of $\omega^{m+1}$.
We usually write ``$\phi_e(x)\converges=y$'' to indicate
that the computation of $\phi_e$ on input $x$ halts
and outputs $y$, and so $\phi_e(x)\converges$ iff $x\in W_e$;
otherwise we write $\phi_e(x)\diverges$.
Also, if the computation halts within $s$ steps,
we write $\phi_{e,s}(x)\converges$.  The set
$W_{e,s}$ is the domain of $\phi_{e,s}$,
so $W_e=\cup_s W_{e,s}$.  Every set $W_{e,s}$
is computable (although the union $W_e$ may not be),
and we take it as a convention of our computations
that only numbers $\leq s$ lie in $W_{e,s}$.

More generally, we define the $\Sigma^0_n$ formulas
by induction on $n$.  The $\Sigma^0_0$ formulas are those
formulas with free variables $x_1,\ldots,x_m$ which define computable
subsets of $\omega^m$.  A $\Pi^0_n$ formula
is the negation of a $\Sigma^0_n$ formula (so
a $\Pi^0_1$ formula is universal, in the
same sense that a $\Sigma^0_1$ formula is existential),
and a $\Sigma^0_{n+1}$ formula in the variable $x$
of the form
$$ \exists y_1\cdots\exists y_m~R(x,y_1,\ldots,y_m),$$
where $R$ is a $\Pi^0_n$ formula.  Thus the subscript
counts the number of quantifier alternations.
(The superscript, often omitted here, refers to the
fact that we quantify only over natural numbers,
not over sets of naturals, or sets of sets of them, etc.)
Consecutive like quantifiers can be collapsed to a single
quantifier, by the use of computable pairing functions:  for each $n$,
there is a computable bijection $\omega^n\to\omega$,
denoted by letting $\la y_1,\ldots,y_n\ra\in\omega$
be the image of the tuple $( y_1\ldots,y_n)\in\omega^n$.
Thus the $\Sigma^0_{n+1}$ formula above could be
expressed as $\exists y R'(x,y)$, where $R'(x,\la y_1,\ldots,y_n\ra)$
holds iff $R(x,y_1,\ldots,y_n)$ holds.  Details can be found in
\cite{Notices08}, where a more general computable
bijection onto $\omega$ from the set $\omega^{<\omega}$ of \emph{all}
finite tuples of natural numbers is also given.  This latter bijection
allows us to use  a single quantifier to quantify over all
polynomials in $F[X]$, for instance, for any computable field $F$
(as defined below).

Turing reducibility and $1$-reducibility are ways of comparing
the complexity of subsets $A,B\subseteq\omega$, both defined
in \cite{S87}.  We write $A\leq_T B$ and $A\leq_1 B$, respectively,
to denote that $A$ is no more complex than $B$ under these relations.
It is well known that, for every $n\in\omega$, there exists a set $S$
which is \emph{$\Sigma^0_{n+1}$-complete}:  $S$ itself is
$\Sigma^0_{n+1}$, and every $\Sigma^0_{n+1}$ set $T$ has $T\leq_1 S$.
Likewise, the complement of $S$ is $\Pi^0_{n+1}$-complete.
This is regarded as an exact assessment of the complexity of $S$;
among other things, it ensures that $S$ is \emph{not}
$\Pi^0_{n+1}$, nor $\Sigma^0_{n}$.  It should be noted
that the class of $\Sigma^0_0$ sets and the class of $\Pi^0_0$
sets coincide:  these are the computable sets, and no set
is $\Sigma^0_0$-complete.  A set which is both $\Sigma^0_n$
and $\Pi^0_n$ is said to be $\Delta^0_n$.
Every $\Delta^0_1$ set is $\Delta^0_0$, but for greater $n$
these classes no longer coincide. The $\Delta^0_{n+1}$
sets are exactly those which are Turing-reducible to a
$\Sigma^0_n$-complete oracle set.  As a canonical
$\Sigma^0_n$-complete set, we usually use
$\emptyset^{(n)}$, the $n$-th jump of the empty set,
as defined in \cite{S87}.

Turing reducibility $\leq_T$ is a partial pre-order
on the power set $\mathcal{P}(\omega)$.  We define $A\equiv_T B$,
saying that $A$ and $B$ are \emph{Turing-equivalent},
if $A\leq_T B$ and $B\leq_T A$.  The equivalence
classes under this relation form the \emph{Turing degrees},
and are partially ordered by $\leq_T$.
In fact, they form an upper semi-lattice under $\leq_T$,
with least element $\bfz$, the degree of the computable sets,
but no greatest element.
One often speaks of a set $A$ as being
\emph{computable in a Turing degree $\bfd$},
meaning that for some (equivalently, for every)
$B\in\bfd$ we have $A\leq_T B$.

\section{Useful Results on Fields}
\label{sec:fields}

Computable fields fit the general definition of
computable structures, which is the basis of computable model theory.

\begin{defn}
\label{defn:computablestructure}
A structure $\mathcal{S}$ in a finite language is \emph{computable} if
its domain is an initial segment of $\omega$, the set of
natural numbers, and all functions and relations in $\mathcal{S}$
are computable when viewed as functions and relations on $\omega$.

A structure $\mathcal{M}$ is \emph{computably presentable}
if it is isomorphic to a computable structure $\mathcal{S}$, in which
case we call $\mathcal{S}$ a \emph{computable presentation} of $\mathcal{M}$.
\end{defn}

A computable field $F$ therefore has domain $\omega$,
or else $\{ 0,1,\ldots,p^k-1\}$, with computable field operations.
Since the symbols $0$ and $1$ have their own meaning
in field theory, we often write the domain of $F$
as $\{ x_0,x_1,\ldots\}$ for clarity.
If the $i$-th partial computable function $\phi_i$ computes
addition on $F$ (so $x_m+x_n=x_{\phi_i(m,n)}$)
and $\phi_j$ computes multiplication on $F$,
then $\la i,j\ra$ is an \emph{index} for $F$.
These definitions are standard in
computable model theory, and we will maintain them here,
but they complicate the discussion of computable fields.
Herewith our conventions.  The standard symbols
$+$, $-$, $\cdot$, $x^n$, and $\frac{x}{y}$,
from field theory all refer to the (computable) operations
in $F$, with ``$-$'' denoting both subtraction and negation as usual.
Likewise, $0$ and $1$ denote the identity elements
of $F$, rather than naming the first two elements
of the domain.  These will be far more useful for us
than the constants and operations on $\omega$ would be.
On the other hand, we use the symbol $<$ to denote
the usual relation on $\omega$, \emph{not on the field $F$}.
Of course, an arbitrary computable field may or may not
be orderable at all, let alone computably orderable,
so field orders (i.e.\ linear orders compatible with
the field operations) will not enter into our discussion.
We will often want to search through the domain $\omega$
until we find an element with a certain property,
and such a search will simply go through the
elements $0,1,2,\ldots$ (or $x_0,x_1,x_2,\ldots$),
using the $<$ relation on $\omega$.
Similarly, phrases such as ``the least element satisfying...''
will mean the least under $<$ on $\omega$.

At a basic level, \cite{Notices08} and \cite{Proceedings08}
are both useful for definitions about computable fields.
They both avoid the notational issue by writing
$\{a_0,a_1,\ldots\}$ in place of $\{0,1,\ldots\}$
as the domain of a field.  For serious research on
computable fields through the twentieth century,
all of \cite{vdW70}, \cite{FS56}, \cite{R60},
\cite{E77}, \cite{MN79}, and \cite{ST99} are
familiar references.

A field is \emph{algebraic} iff it is an algebraic extension
of its prime subfield.  (The prime subfield is just the smallest subfield;
it is a copy of either the rationals $\Q$ or the $p$-element field
$\F_p$, depending on characteristic.)  In this paper, we will restrict
ourselves entirely to algebraic fields, although in characteristic
$0$, our results carry over to the case of fields of finite
transcendence degree over $\Q$, just by fixing a transcendence basis $B$
and treating $\Q(B)$ as the prime subfield, over which the
rest of the field is algebraic.

We restrict ourselves further to the case of algebraic fields with
splitting algorithms.  A computable field $F$ has a
\emph{splitting algorithm} if its \emph{splitting set}
$$ S_F = \set{p\in F[X]}{p\text{~is reducible in~}F[X]}$$
is computable.  (To clarify:  $p$ is \emph{reducible} if it can be expressed
as a product of nonconstant polynomials in $F[X]$;
it need not split into linear factors, but it must split
into at least two proper factors.)
The polynomial ring $F[X]$ may be presented effectively
(i.e.\ as a computable ring) by use of the computable bijection
from $\omega^{<\omega}$ onto $\omega$ described in Section \ref{sec:computability}.
The next result relates $S_F$ to other properties of
the field $F$.  This is a direct consequence of Rabin's Theorem,
first published in \cite{R60} in 1960, and is discussed in more
detail in \cite{M09b}, where it appears as Corollaries 2.7 and 2.8.
\begin{lemma}
\label{lemma:invariant}
Let $F$ be any computable field.  Then the following are all
Turing-equivalent:  the splitting set
of $F$
$$ S_F=\set{p(X)\in F[X]}{p(X)\text{~has a proper factorization in~}F[X]},$$
the \emph{root set} of $F$
$$ R_F=\set{p(X)\in F[X]}{(\exists r\in F) p(r)=0},$$
the \emph{root function} $g_F$ of $F$
$$ g_F(p) = |\set{r\in F}{p(r)=0}|,\text{~with~}\dom{g_F}=F[X],$$
and the \emph{root multiplicity function} of $F$ (which is the same as $g_F$,
except that roots are counted by multiplicity).
Moreover, the Turing reductions are uniform in (an index for) $F$.

Moreover, any two isomorphic computable algebraic
fields $F$ and $\Ftilde$ must have
Turing-equivalent splitting sets, and the Turing reductions
are uniform in $F$ and $\Ftilde$.
Hence, for computable algebraic fields, the Turing degree of each item
above is an invariant of the isomorphism type.
\qed\end{lemma}

Notice that an algebraic field
need not be finitely generated (equivalently, need not have finite degree)
over the prime subfield.  Indeed, finitely generated
computable fields are relatively straightforward objects:
for one thing, they always have splitting algorithms,
as do all prime fields.  This was shown in 1882 by Kronecker,
\label{Kronecker}
in \cite{K1882}; a discussion in modern terms appears
in \cite[\S 2]{M09b}, where it is also explained
how we can determine splitting algorithms, uniformly
in the generators, for all finitely generated subfields $E$
of any computable algebraic field $F$, and how we can use these
to determine the Galois group of any finitely generated subfield
$E\subseteq F$ over any subfield of $E$, uniformly
in the generators of $E$ and the subfield.
We view the Galois group as a permutation group on
the generators and their conjugates over the ground field;
this is a useful way to consider automorphisms of $E$
as finitary objects rather than as functions.

We add a new notion, the conjugacy function, which is related but not always
Turing-equivalent to the sets and functions from Lemma \ref{lemma:invariant}.
\begin{defn}
\label{defn:conjfct}
Let $F$ be a computable field, with prime subfield $Q$.
The \emph{conjugacy function for $F$} is the function
$h:F\to\omega$ defined by:
$$ h(x) =\left\{\begin{array}{l}
\diverges,\text{~~~~if $x$ is transcendental over~}Q\\
|\set{y\in F}{x\text{~and $y$ are conjugate over~}Q}|,\text{~else.}
\end{array}\right.$$
\end{defn}
In general, the conjugacy function is a partial function.  Normally
its initial step is to find the minimal polynomial of $x$ over $Q$,
and if $x$ is transcendental, then this search will never converge.
In practice we are usually concerned with the case of an
algebraic extension $F$ over $Q$, in which case $h$ is total.

The conjugacy function is always computable from the splitting set,
or from any of the other Turing-equivalent sets in Lemma \ref{lemma:invariant}.
The splitting set allows us to find the minimal polynomial $q(X)$ of
any $x$ algebraic over $Q$, and then we use the root function to determine $h(x)$.
However, the conjugacy function may have strictly smaller
Turing degree than those sets.  For example, in any normal
extension of $Q$, the conjugacy function is computable,
being given by the degree of the minimal polynomial;
yet the computable normal algebraic field
$\Q[\sqrt{p_n}~|~n\in K]$, generated by the square roots of
the primes $p_n$ with $n$ in the halting set $K$, has
splitting set of degree $\bfz'$.

The analogue of Rabin's Theorem for the conjugacy function
is the following.
\begin{prop}
\label{prop:Rabinconj}
Let $F$ be any computable algebraic field with prime subfield $Q$,
and $g$ any computable field embedding of $F$ into a computable
presentation $\Qbar$ of the algebraic closure of $Q$.
Then the conjugacy function $h$ of $F$ is computable
iff the image $g(F)$ is a computable subfield of its
normal closure within $\Qbar$, i.e.\ iff
there exists a partial computable function $\psi$
whose domain is the normal closure of $g(F)$
within $\Qbar$, such that $\psi$ is the characteristic
function of $g(F)$ on that domain.
More generally, this $\psi$ is Turing-equivalent to $h$.
\end{prop}
Notice that the normal closure of $g(F)$ within $\Qbar$
is independent of the choice of the computable embedding $g$.
\begin{pf}
Suppose first that we have an $h$-oracle.  Given any $x\in\Qbar$,
$\psi$ waits until some $y$ appears in $F$ such that $g(y)$
and $x$ have the same minimal polynomial over the prime
subfield $Q$ of $\Qbar$.  If this never happens, then $x$
does not lie in the normal closure of $g(F)$, so $\psi(x)$
need not converge.  If $y$ appears, then we compute $h(y)$,
using the oracle, and find all conjugates
$y=y_1,y_2,\ldots,y_{h(y)}$ of $y$ over $Q$ in $F$.
Then $x\in g(F)$ iff $(\exists i\leq h(y))~g(y_i)=x$.

Conversely, with an oracle for $\psi$, the program for $h$
accepts any $y\in F$ as input, computes
the minimal polynomial $q(X)$ of $g(y)$ over $Q$,
and finds all roots $x_1=g(y),x_2,\ldots,x_{\deg{q}}$ of
$g(y)$ in $\Qbar$.  Then $\psi(x_i)\converges$ for all $i\leq\deg{q}$,
and $h(y)$ is the number of these $i$ for which $\psi(x_i)=1$.
\qed\end{pf}

Of course, often one wants to know about conjugates over ground fields
other than the prime subfield.  If this ground field is algebraic,
then the conjugacy function computes this information.
\begin{lemma}
\label{lemma:genconj}
Let $F$ be a computable field, with prime subfield $Q$
and conjugacy function $h$.  Then for any subfield $E\subseteq F$
algebraic over $Q$, the function
$$ h_E(x) =\left\{\begin{array}{l}
\diverges,\text{~~~~if $x$ is transcendental over~}E\\
|\set{y\in F}{x\text{~and $y$ are conjugate over~}E}|,\text{~else}
\end{array}\right.$$
is computable uniformly in oracles for $h$ and the splitting set of $E$.
In particular, if $E$ is generated by a finite set $z_1,\ldots,z_n$ over $Q$,
then we can compute $h_E$ uniformly in $h$ and the tuple $\zvec$.
\end{lemma}
\begin{pf}
For any $x\in F$, the splitting set for $E$ allows us to find the minimal polynomial
$p(X)$ of $x$ over $E$.  We also determine the minimal polynomial $q(X)$
of $x$ over $Q$, and search until we have found all roots $x=x_1,x_2,\ldots,x_{h(x)}$
of $q(X)$ in $F$.  Now $p(X)$ must divide $q(X)$ in the ring $E[X]$, so
$h_E(x)$ is just the number of these roots $x_i$ satisfying $p(x_i)=0$.

Here $F$ is not assumed to be algebraic over $Q$, but $E$ is algebraic.
So, if $x\in F$ is transcendental over $Q$, then $x$ is also transcendental
over $E$, and therefore the computation of $h_E(x)$ described above
never halts, which is exactly the prescribed outcome.
\qed\end{pf}

Next we introduce several standard facts about fields
which do not involve computability.
The following result appears as Lemma 2.10 in \cite{M09b}.
\begin{lemma}
\label{lemma:endo}
For an algebraic field $F$,
every endomorphism (i.e.\ every homomorphism from $F$ into
itself) is an automorphism.
\qed\end{lemma}
\begin{cor}
\label{cor:endo}
If $F\cong\Ftilde$ are isomorphic algebraic fields, and $f:F\to\Ftilde$
is a field embedding, then the image of $f$ is all of $\Ftilde$.
That is, such an $f$ must be an isomorphism.
\qed\end{cor}

Before continuing to Proposition \ref{prop:Konig}, we explain
the intuition behind the \emph{embedding tree},
which will be stated formally in Definition \ref{defn:isotree},
but is used here as well.  To construct a field embedding
from $E$ into $F$ effectively, one usually begins with
the prime subfield $Q$ of $E$ (which is always computably enumerable
within $E$, and can be mapped effectively onto the
prime subfield of $F$, provided the fields have the same
characteristic).  Then one extends this partial embedding $f_0$
to the least element $x_0$ in the domain $\{ x_0,x_1,\ldots\}$ of $E$
(which might already lie in $Q$, of course), then the second-least,
and so on.  The options for the image of $x_{s+1}$ depend
on the choices we made for the images of $x_0,\ldots,x_s$,
naturally.  At each stage,
the number of ways to extend the embedding
$f_s:Q(x_0,\ldots,x_{s-1})\embeds F$ to $x_s$
is bounded by the degree of the minimal polynomial
of $x_s$ over $Q(x_0,\ldots,x_{s-1})$, and there
might be no way at all to do so, even assuming
that $E$ embeds into $F$, because we may
have made bad choices at previous stages.
Our choices thus naturally give rise to the
\emph{embedding tree}, with a root node
representing the partial embedding $f_0$
(to which every embedding of $E$ into $F$ restricts),
and with each node $\sigma$ of length $|\sigma|=s$
having one immediate successor for each of the
finitely many possible images in $F$ of $x_s$,
given our previous choices.  Hence one may view $\sigma$
as a node in the tree $F^{<\omega}$ of finite
sequences of elements of $F$:
the sequence $\la y_0,\ldots,y_{s-1}\ra$
describes the map with $x_i\mapsto y_i$ for all $i<s$,
and this sequence lies in the isomorphism tree
iff that map extends to an embedding of $Q(x_0,\ldots,x_s)$
into $F$.  The infinite paths
through this tree correspond precisely to the embeddings
of $E$ into $F$.  Moreover, the tree itself is computable,
in the sense that we can decide exactly which
sequences in $F^{<\omega}$ lie in the tree.
Nevertheless, some or all of the paths through
the tree may be noncomputable.

\begin{prop}
\label{prop:Konig}
Let $E$ and $F$ be algebraic field extensions of a
common subfield $Q$.  Then $E$
embeds into $F$ over $Q$ iff every finitely generated subfield of $E$
containing $Q$ embeds into $F$ over $Q$.
\end{prop}
\begin{pf}
The key to the nontrivial direction is K\"onig's Lemma,
applied to the \emph{embedding tree} $T_{E,F; Q}$ for $E$ and $F$
over $Q$, as we now explain.  Write the domain of $E$ as
$E=\{ x_0,x_1,\ldots\}$, with $\omega$ as the domain of $F$.
For each $s$, pick any $q_s\in Q[X_0,\ldots,X_s]$
such that $q_s(x_0,\ldots,x_{s-1},X)$ is the minimal
polynomial of $x_s$ over $Q[x_0,\ldots,x_{s-1}]$.
We consider the set $T_{E,F;Q}$ of finite sequences $\sigma$  of
natural numbers:
$$ T_{E,F;Q} =\set{\sigma}{(\forall s<|\sigma|)~q_s(\sigma(0),\ldots,\sigma(s))=0}.$$
Notice that if $\sigma=(a_1,\ldots,a_n)\in T_{E,F;Q}$, then so is
$\tau=(a_1,\ldots,a_m)$ for each $m=|\tau|\leq n=|\sigma|$.
That is, $T_{E,F;Q}$ is closed under initial segments, so we
view it as a tree.

Sequences $\sigma\in T_{E,F;Q}$ correspond to embeddings
of $Q[x_0,\ldots,x_{|\sigma|-1}]$ into $F$,
since for each $s$ less than the length $|\sigma|$ of the sequence,
$\sigma(s)$ satisfies the same minimal polynomial
over $Q(\sigma(0),\ldots,\sigma(s-1))$ in $F$
that $x_s$ satisfies over $Q(x_0,\ldots,x_{s-1})$ in $E$.
Clearly $T_{E,F;Q}$ is finite-branching:
the number of successors of $\sigma$ is at most the degree
of the variable $X_{|\sigma|}$ in $q_{|\sigma|}$.
(A more general formal description appears as Definition
\ref{defn:isotree} below.)

K\"onig's Lemma states that every finite-branching tree
with infinitely many nodes must contain an infinite path.
By assumption, every subfield $Q[x_0,\ldots,x_n]$ embeds into $F$,
so $T_{E,F;Q}$ contains a node
of each length $n\in\omega$.  Therefore, $T_{E,F;Q}$
contains an infinite path,
which defines an embedding of $F$ into $E$.
\qed\end{pf}

A less slick but intuitively clearer proof of this result
(with $Q$ as the prime subfield) appears
within the proof of \cite[Appendix A, Thm.\ 2]{FKM09}.
Sometimes we will write just $T(E,F)$, in which case $Q$
is to be understood as the prime subfield of both.
The concept of the embedding tree foreshadows the
\emph{isomorphism trees} in Definition \ref{defn:isotree} below.
\begin{cor}
\label{cor:Konig}
Two algebraic fields $E$ and $F$ are isomorphic
over a common subfield $Q$ iff every
finitely generated subfield (containing $Q$) of each one
embeds over $Q$ into the other.
\end{cor}
\begin{pf}
By Proposition \ref{prop:Konig}, the latter condition is equivalent
to $E$ and $F$ both embedding into each other over $Q$.  But then the
composition of these two embeddings is an automorphism,
by Lemma \ref{lemma:endo}.
\qed\end{pf}
The upshot of this corollary is that often, to build a computable
field $F$ isomorphic to a given computable field $E$,
we can simply construct $E$ as a union of nested,
uniformly computable fields $E_s\subseteq E_{s+1}$
with each $E_s$ isomorphic to the subfield $F_s\subseteq F$
generated by the first $s$ elements of $F$.  There is no need for
the isomorphisms $f_s:F_s\to E_s$ to have a limit;
Corollary \ref{cor:Konig} does all the work for us.

For this paper we have two new definitions, arising out of the
standard concept of conjugacy.  Examples appear in
Section \ref{sec:EX}.
\begin{defn}
\label{defn:conjugate}
Let $E\subseteq F$ be any field extension.  Two elements
$a,b\in F$ are \emph{conjugate over $E$} if they have the
same minimal polynomial in $E[X]$.  It is well known that then
the subfields $E[a]$ and $E[b]$ are isomorphic, via a map
fixing $E$ pointwise with $a\mapsto b$.  We say that $a$ and $b$ are
\emph{true conjugates in $F$ over $E$} if there exists
an automorphism $\psi$ of $F$ with $\psi(a)=b$ and
$\psi\res E$ being the identity map.  If $a$ and $b$
are conjugate over $E$ but are not true conjugates
in $F$ over $E$, we call them \emph{false conjugates
in $F$ over $E$}.
\end{defn}

Finally, we give the computability-theoretic version
of the classical Theorem of the Primitive Element.
\begin{thm}[Effective Theorem of the Primitive Element]
\label{thm:primitive}
Let $E$ and $F$ be computable fields, with $E\subseteq F$
a separable algebraic extension and with $E$ c.e.\ as a subfield
of $F$.  Then for any elements $x_1,\ldots,x_n\in F$,
we can effectively find a \emph{primitive generator} $y\in F$
for these elements.  That is, we can find a $y\in F$
such that $E[y]=E[x_1,\ldots,x_n]$.  The procedure
for finding $y$ is uniform in $n$,
the generators $\la x_1\ldots,x_n\ra$, the enumeration of $E$
within $F$, and the field operations in $F$.
\end{thm}
\begin{pf}
The existence of such an element $y$ is the classical theorem
see e.g.\ \cite[p.\ 139]{vdW70}, and is made effective
by a direct search for $y$.  An arbitrary $y\in F$ generates $E[x_1,\ldots,x_n]$ iff
$$(\exists p\in E[X_1,\ldots,X_n])
(\exists q_1,\ldots,q_n\in E[Y])[y=p(x_1,\ldots,x_n)~\&~
\text{all~}x_i=q_i(y)],$$
and so we can find $y$, with the uniformities described above.
For a proof giving an actual formula for the generator $y$,
see \cite{FJ86}.
\qed\end{pf}

\section{A First Example}
\label{sec:EX}

We start with an example, to introduce the concepts
that will be used later in our analysis of computable
categoricity for algebraic fields.
Let $F_0$ be a computable presentation of
the normal algebraic extension of the rationals by the square roots
of all rational primes:  $F_0=\Q[\sqrt{p_0},\sqrt{p_1},\ldots ]$.
The domain of $F_0$ is $\omega$, and we use $\sqrt{p_i}$
to refer to the lesser of the two square roots of $p_i$
in $F_0$, under the ordering $<$ on the domain $\omega$.
Now, for any $W\subseteq\omega$, let $F_W$ be the extension
of $F_0$ in which we adjoin a square root of $\sqrt{p_i}$
iff $i\in W$.  Notice that no $-\sqrt{p_i}$ acquires a square root
of its own during this process, whether or not $i\in W$.
(To see this, embed $F_W$ into the field $\mathbb{R}$.)
So $F_W$ is not normal over $\Q$ (unless $W=\emptyset$),
although it is normal over $F_0$, which in turn is normal over $\Q$.
$F_W$ is computably presentable
iff $W$ is computably enumerable,
in which case we take $F_W$ to denote a computable presentation
built over our original presentation of $F_0$.  The domain of $F_0$
can be the set of even elements of $\omega$, for instance,
with the odd elements added as numbers $i$ appear in $W$
and dictate the adjoinment of square roots of $\sqrt{p_i}$
and the new elements they generate.

Notice that for all $i\in W$, $\sqrt{p_i}$ and $-\sqrt{p_i}$
are false conjugates in $F_W$ over $\Q$, as in
Definition \ref{defn:conjugate}:  they both have minimal
polynomial $X^2-p_i\in \Q[X]$, but neither can be mapped to the
other by any automorphism of $F_W$ over $\Q$, since
$\sqrt{p_i}$ has a square root of its own and $-\sqrt{p_i}$
does not.  On the other hand, for $i\notin W$, $\pm\sqrt{p_i}$
are true conjugates in $F_W$ over $\Q$, since there is an
automorphism of $F_W$ mapping one to the other.

Before stating the main result for these fields $F_W$,
we recall the concept of computable inseparability.
\begin{defn}
\label{defn:insep}
Two sets $P$ and $N$ are \emph{computably inseparable}
if there is no computable set $C$ with $P\subseteq C$
and $N\subseteq\overline{C}$.
\end{defn}
A standard example (see \cite[I.4.22]{S87}) has
$P=\set{n\in\omega}{\phi_n(n)\converges =0}$ and
$N=\set{n\in\omega}{\phi_n(n)\converges =1}$.
(Recall that ``$\phi_e(x)\converges=0$'' means that the $e$-th
partial computable function, when run with the input $x$,
halts and outputs $0$.)
The following result
was first proven by Yates, who saw that it followed from a
construction of Friedberg; but the proof was only published
by Cleave \cite{C70} in 1970, some years later.
\begin{thm}[Friedberg \& Yates \cite{C70}]
\label{thm:Yates}
Every noncomputable c.e.\ set is the union
of two disjoint, computably inseparable c.e.\ subsets.
\qed\end{thm}
Conversely, if $W$ is computable, then
for every partition of $W$ into c.e.\ sets $A$ and $B$,
both $A$ and $B$ are computable, since the complement $\overline{A}$
of the c.e.\ set $A$ is the set $(\overline{W}\cup B)$, which is also c.e.

\begin{prop}
\label{prop:EX}
For any c.e.\ $W$, the computable field $F_W$
is computably categorical iff $W$ cannot be partitioned
into computably inseparable c.e.\ subsets.
\end{prop}
\begin{cor}
\label{cor:EX}
$F_W$ is computably categorical iff $W$ is computable.
\qed\end{cor}
The corollary is immediate, using Theorem \ref{thm:Yates}.
The point of Proposition \ref{prop:EX} is its proof (below),
rather than the result itself,
since the proof illustrates the usefulness of the concept
of computable inseparability (see Definition \ref{defn:insep}).
\begin{pf}
First suppose that there is no partition of $W$ into
computably inseparable c.e.\ sets, and let
$\Ftilde$ be a computable field isomorphic to $F_W$.
Define
\begin{align*}
P&=\set{i\in\omega}{\exists x,y,z\in \Ftilde
[x<y~\&~x^2=y^2=\ptilde_i~\&~z^2=x]}\\
N&=\set{i\in\omega}{\exists x,y,z\in \Ftilde
[y<x~\&~x^2=y^2=\ptilde_i~\&~z^2=x]}.
\end{align*}
(Here $<$ refers to the standard order on the domain
$\omega$ of $\Ftilde$, but $\ptilde_i$ is the $i$-th rational prime
in $\Ftilde$, and squares refer to the field multiplication in $\Ftilde$.)
The computability of $\Ftilde$ shows that $P$ and $N$ are both c.e.,
and since $\Ftilde\cong F_W$, we know that $P\cup N=W$
and $P\cap N=\emptyset$.
By assumption, therefore, there exists a computable $C$
with $P\subseteq C$ and $N\subseteq\overline{C}$.
We define our computable isomorphism $f:F_W\to\Ftilde$
beginning with the computable subfield $F_0\subseteq F_W$.
$f$ is uniquely determined on $\Q$ within $F_0$.
If $i\in C$, then we map
$\sqrt{p_i}$ (from $F_W$) to the lesser square root
of $\ptilde_i$ in $\Ftilde$, with $-\sqrt{p_i}$ therefore
mapped to the greater one.  In this situation we know that
either $i\notin W$ (in which case $\sqrt{p_i}$ can be mapped
to either square root of $\ptilde_i$) or $i\in P$
(in which case we made the correct choice by mapping $\sqrt{p_i}$
to the lesser square root of $\ptilde_i$).
If $i\notin C$, we do the opposite,
mapping $\sqrt{p_i}$ to the greater square root
of $\ptilde_i$ in $\Ftilde$.  With $i\notin C$,
we know that either $i\notin W$ or $i\in N$,
and so again the choice we made was a correct choice.
This much is readily computable, since $C$ is computable.

Moreover, as elements $x$ outside of $F_0$ appear
in $F_W$, we can compute the (unique) extension of $f$
to those elements:  if a square root $z_i$ of $\sqrt{p_i}$
ever appears in $F_W$, then $i\in W=P\cup N$, and our choice
of $f(\sqrt{p_i})$ using $C$ ensures that $f(\sqrt{p_i})$
has a square root of its own in $\Ftilde$, to which we map $z_i$.
Since $x$ must have been generated by finitely many of these
$z_i$, this will eventually allow us to extend $f$ to $x$.
Clearly, then, $f$ is a computable embedding of $F_W$ into $\Ftilde$,
and by Corollary \ref{cor:endo}, $f$ must be an isomorphism.
Thus $F_W$ is computably categorical.

Now suppose there exist disjoint, computably inseparable
sets $P$ and $N$ whose union equals $W$.  We define a computable
field $F$ isomorphic to $F_W$ by starting with $F_0$ and enumerating
$P$ and $N$.  Whenever any $i$ appears in $P$, we adjoin a
square root of $\sqrt{p_i}$ to $F$, while when an $i$ enters $N$,
we adjoin a square root of $-\sqrt{p_i}$ to $F$.  Since $P$ and
$N$ form a partition of $W$, this $F$ must be isomorphic to $F_W$.
However, if $g:F_W\to F$ is any isomorphism, then
$C=\set{i\in\omega}{g(\sqrt{p_i}) < g(-\sqrt{p_i})}$
is computable in $g$ and must contain $P$ while not
intersecting $N$.  (Here again $<$ is the standard order
on the domain $\omega$ of $F$.)
Therefore no such $g$ can be computable, and $F$
is a computable field isomorphic to $F_W$ but not computably
isomorphic to it.
\qed\end{pf}
Relativizing this proof yields an immediate corollary.
\begin{cor}
\label{cor:catspec}
For c.e.\ sets $W$, the above field $F_W$
is $\bfd$-computably categorical iff
the degree $\bfd$ computes separations
of every partition of $W$ into two c.e.\ subsets.
That is, the degrees with this property form the
\emph{categoricity spectrum} of $F_W$,
as defined in \cite{FKM08}.
\qed\end{cor}
We note that the splitting set of this field $F_W$
is Turing-equivalent to $W$, hence is in general noncomputable.
There is a similar procedure which builds a computable field $K_W$
with computable splitting set, and with the same categoricity spectrum:
the basic strategy is that, if
$n$ enters $W$ at stage $s$, one should adjoin a root
of a polynomial $q_n(\sqrt{p_n},Y)$ in such a way
that $q_n(-\sqrt{p_n},Y)$ has no root.  By ensuring
that these polynomials all have distinct prime
degrees, we can make certain that $q_n(-\sqrt{p_n},Y)$
never acquires a root at any other stage,
and by making that prime degree $d$ be $>s$, we can
keep the splitting set of $F$ computable.
Unfortunately, $q_n(\sqrt{p_n},Y)$ cannot simply
say $Y^d=\sqrt{p_n}$, because then the negative
of this root would be a $d$-th root of $-\sqrt{p_n}$.
For the existence of the requisite polynomials $q_n$,
see \cite[Prop.\ 2.15]{M09b}, and for the full construction
see Theorem 3.4 there.  With this note, we have informally
established the following corollary.

\begin{cor}
\label{cor:Sigma03}
Computable categoricity for computable algebraic fields
with splitting algorithms is $\Sigma^0_3$-hard.
\end{cor}
\begin{pf}
This means that there is a $1$-reduction $f$ from some
$\Sigma^0_3$-complete set $S$:  for all $e\in S$, $f(e)$
is the index of an algebraic field with a splitting algorithm,
which is computably categorical iff $e\in S$.  We choose
the set $\textbf{Rec}=\set{e}{W_e\leq_T\emptyset}$ to serve
as $S$, and $f(e)$ is the field index produced by the above
construction on $W_e$.  For details, see \cite[\S IV.3]{S87}.
\qed\end{pf}

It was also shown in \cite{M09b} that a field such as $K_W$
(with $W>_T\emptyset$),
with a splitting algorithm but not computably categorical,
cannot have any least Turing degree in its categoricity
spectrum.  By Corollary \ref{cor:catspec}, this shows
that for noncomputable c.e.\ sets $W$, there is no least degree
which computes a separation of every partition of $W$
into two c.e.\ subsets.  Of course, this result can be shown
directly, without the excursus into computable model theory!
The key to the direct proof, and also to the result in \cite{M09b},
is the Low Basis Theorem of Jockusch and Soare,
from \cite{JS72}.

%

\section{Isomorphism Trees for Fields}
\label{sec:isotrees}

Computable separability arises frequently in discussions
of computable trees and paths through them.  In light of the connection
seen in the preceding section between computable separability
and computably categorical fields, it is not surprising that,
in discussing isomorphisms between fields, we will make
great use of the \emph{isomorphism tree}.  This concept
was introduced in \cite[\S 5]{M09b} and used in
Proposition \ref{prop:Konig} above; it may be useful
for the reader to refer back now to the discussion there.
Here we give a full definition, generalizing the work in \cite{M09b}
to the case where the two fields are not necessarily isomorphic.
$F^{<\omega}$ denotes the set of finite sequences $\sigma$ of
elements of $F$ (or, formally, functions $\sigma:\{0,1,\ldots,|\sigma|-1\}\to F$,
where $|\sigma|$ is the length of $\sigma$).
\begin{defn}
\label{defn:isotree}
Let $E$ and $F$ be computable algebraic fields of the same characteristic,
and let $\la x_0,x_1,\ldots\ra_{i\in J}$ be a computable sequence
(finite or infinite) of elements of $E$ which together generate $E$
over its prime subfield $Q$.
Let $\la q_i\ra_{i\in J}$ be polynomials such that every
$q_i\in Q[X_0,\ldots,X_i]$ and for each $i$, $q_i(x_0,\ldots,x_{i-1},X)$
is the minimal polynomial of $x_i$ over $Q(x_0,\ldots,x_{i-1})$.
The \emph{embedding tree} $I_{EF}$ is the computable tree
$$ I_{EF}=\set{\sigma\in F^{<\omega}}{(\forall i<|\sigma|)
\qtilde_{i}(\sigma(0),\ldots,\sigma(i))=0},$$
where $\qtilde_i$ is the image of $q_i$ under the unique
embedding of $Q$ into $F$.
The height of $I_{EF}$
is $|J|+1$ for finite generating sets $J$, or $\omega$
if $|J|=\omega$, and $I_{EF}$ must be finite-branching:
for each $\sigma\in I_{EF}$ of any length $n$, there are only finitely many
$\tau\in I_{EF}$ such that $\sigma\subseteq\tau$
and $|\tau|=n+1$, because there are only finitely
many roots of $q_{n+1}(\sigma(0),\ldots,\sigma(n-1),X)$ in $F$.

When $E$ and $F$ are known to be isomorphic,
the term \emph{isomorphism tree} is often used
(with Corollary \ref{cor:endo} as justification).
If $E=F$ and $I=\la x_0,x_1,\ldots\ra$ simply lists the domain of $F$ in order
(as elements of $\omega$), then we write $I_F$ for $I_{FF}$ and call this the
\emph{automorphism tree} of $F$.  (We sometimes abuse this terminology
by writing $I_F$ even when a different generating sequence is
being used, as long as it is clear which sequence it is.)
The \emph{identity path} in the automorphism
tree $I_F$ contains all nodes of the form $\la x_0,x_1,\ldots,x_i\ra$;
clearly this is a path through $I_F$, corresponding to the identity
automorphism, since the map it defines sends each $x_i$ to $x_i$ itself.
The \emph{level function}
for $F$ is the function giving (for each $n$) the number
of nodes at level $n$ in $I_F$.
\end{defn}

In Definition 5.1 in \cite{M09b}, the sequence $\la x_i\ra_{i\in I}$
was assumed simply to be the domain of $F$.  For the purposes
of this paper, we generalized to arbitrary computable generating sequences,
but this does not change the central point about isomorphism trees,
which is the following.
\begin{lemma}
\label{lemma:paths}
If $F\cong\Ftilde$, then there is a bijection $\psi\mapsto h_\psi$ from the set of
field isomorphisms $\psi:F\to\Ftilde$ onto the set of paths through $I_{F\Ftilde}$,
with $\psi\equiv_T h_\psi$ for all such isomorphisms $\psi$.
(For each path $h\in\Ftilde^\omega$ through $I_{F\Ftilde}$,
we write $\psi_h$ for the unique field
isomorphism such that $h=h_{\psi_h}$, namely $\psi_h(x_i)=h(i)$.)
\end{lemma}
\begin{pf}
For each $\psi$, the path $h_\psi$ is the set of all nodes $\sigma\in I_{F\Ftilde}$
of the form $\la \psi(x_0),\psi(x_1),\ldots, \psi(x_{|\sigma|-1})\ra.$
Since $\psi$ is an isomorphism, this is clearly a path in $I_{F\Ftilde}$,
and conversely, for any path $h$, setting $\psi(x_i)=\sigma(i)$
for the unique $\sigma\in h\cap\Ftilde^{(i+1)}$ defines a field embedding $\psi$,
which must be an isomorphism by Corollary \ref{cor:endo},
with $h=h_\psi$.  The Turing equivalence is immediate
from these definitions.  Also, when $\Ftilde=F$, the
image $\hid$ of the identity automorphism
is the identity path through $I_F$.
\qed\end{pf}

We note that Lemma \ref{lemma:paths} holds
independently of the choice of generating sequence of $F$.
The isomorphism trees for two different generating sequences
may not be isomorphic to each other as trees,
but nevertheless the set of paths through each one
still corresponds to the set of isomorphisms from $F$ onto $\Ftilde$.
The polynomials $q_i$ are not uniquely determined by
the generating sequence; one can use any $q_i$,
satisfying the condition that $q_i(x_0,\ldots,x_{i-1},X)$
be the minimal polynomial of $x_i$ over $Q(x_0,\ldots,x_{i-1})$,
and $I_{F\Ftilde}$ will come out exactly the same.
Finally, replacing $\Ftilde$ by a different field
isomorphic to $F$ would not change anything about
the isomorphism tree except the names of the nodes.
\begin{lemma}
\label{lemma:autotree}
For computable algebraic fields $F\cong\Ftilde$, and for a fixed
generating sequence for $F$, we always have an isomorphism of trees:
$$ I_{F\Ftilde}\cong I_F.$$
Indeed, there is a canonical 1-1 map from field isomorphisms
$\psi:F\to\Ftilde$ to tree isomorphisms $H_\psi: I_F\to I_{F\Ftilde}$,
i.e.\ to bijective maps $H$ which preserve
the successor relation on the trees.
Hence there is also a canonical 1-1 map from paths through $I_{F\Ftilde}$ to such tree isomorphisms,
with $H_\psi\equiv_T\psi$ uniformly in $\psi$.
\end{lemma}
\begin{pf}
Fix any field isomorphism $\psi$, and, for $\sigma\in I_F$
with $|\sigma|=n$, define $H_\psi$ on $I_F$ by
$$ H_\psi (\sigma) = \la \psi(\sigma(0)),\ldots,\psi(\sigma(n-1))\ra\in I_{F\Ftilde}.$$
Clearly $H_\psi$ is a tree isomorphism: $H_\psi(\sigma)\in I_{F\Ftilde}$ because
$$ q_n( \psi(\sigma(0)),\ldots,\psi(\sigma(n-1)))=
\psi (q_n(\sigma(0),\ldots,\sigma(n-1))) =\psi(0)=0,$$
using the definition of $I_F$ to see that $q_n(\sigma(0),\ldots,\sigma(n-1))=0$,
and bijectivity and preservation of the successor relations on the trees
follow from $\psi$ being a field isomorphism.

The definition of $H_\psi$ shows it to be computable from $\psi$.
Conversely, if we know $H_\psi$, we can apply it to nodes
$\la x_0,\ldots,x_s\ra$ on the identity path in $I_F$,
yielding $\la \psi(x_0),\ldots,\psi(x_s)\ra$, thus
computing $\psi$ on the generating sequence $\la x_s\ra_{s\in I}$
for the field $F$.
\qed\end{pf}

The canonical extension of the pairing need not be the only possible
extension to an isomorphism from $I_F$ onto $I_{F\Ftilde}$,
and so the 1-1 map in Lemma \ref{lemma:autotree} need
not be onto.  Any tree isomorphism $H:I_F\to I_{F\Ftilde}$
must map the identity path to some path $h$ through $I_{F\Ftilde}$,
and thus gives rise to a field isomorphism $\psi_h$, although
$H_{\psi_h}$ may not equal $H$.  For example,
if $F$ has generating sequence $x_0,x_1, x_2$,
where $x_0\in\R$ and $x_2\notin\R$
are cube roots of $2$ and $x_1$ is a root of $p(x_0,Y)$
for some $p\in\Q[X,Y]$ such that the resulting field
$F$ contains no roots of $p(x_2,Y)$ or $p(\overline{x_2},Y)$,
then the only nontrivial automorphism of $F$ is complex conjugation,
with $x_2\mapsto\overline{x_2}$, yet the tree $I_F$
for this generating sequence has four automorphisms.
The automorphism tree $I_F$ contains two
terminal nodes $\la x_2\ra$
and $\la\overline{x_2}\ra$ at level $1$,
which represent what happens when one maps $x_0$ to $x_2$
(or to $\overline{x_2}$) and then finds that there is no
element of $F$ to which to map $x_1$.
There is an automorphism $\theta$ of $I_F$ which fixes each path through $I_F$,
but interchanges these two terminal nodes.
This $\theta$ is not equal to $H_\psi$, no matter
which of the two automorphisms of $F$ we try:
$H_{\text{id}}$ fixes $I_F$ pointwise, and when
$\psi$ is complex conjugation, $H_{\psi}$ interchanges
the two paths through $I_F$.  So the tree
isomorphism $\theta$ is not the image of any field isomorphism
under the canonical 1-1 map from Lemma \ref{lemma:autotree}.

The Turing degree of the branching of an automorphism
tree $I_F$ is of interest, especially in light of its use in
\cite{M09b}.  Recall that the level function
for $F$ gives the number of nodes at each level $n$ in $I_F$.
\begin{lemma}
\label{lemma:isotreewidth}
For any algebraic field $F$, the level function $l$ for $F$
is Turing-computable in the conjugacy function $h$ for $F$.
\end{lemma}
\begin{pf}
The number $l(\sigma)$ of immediate successors of a node $\sigma\in I_F$
is exactly the number of conjugates of $x_{|\sigma|}$ in $F$
over the subfield $Q[x_0,\ldots,x_{|\sigma|-1}]$.
We then appeal to Lemma \ref{lemma:genconj}.
\qed\end{pf}

The interesting part, however, is that $l$ need not
be Turing-equivalent to $h$, and indeed the Turing degree
of $l$ need not even be invariant under distinct presentations
of the field $F$.  Thus, Lemma \ref{lemma:invariant}
applies to the root set, the splitting set, the root function,
and the conjugacy function, but \emph{not} to the level function!
To see this, consider the following micro-example.
Suppose that $F$ has domain $\{z_0,z_1,\ldots\}$,
and that $z_0^2=z_5^2=2$, $z_1^2=z_2^2=z_0$,
$z_3^2=z_4^2=z_1$, and that if a certain $m$ appears in the
halting set $K$, then $z_5$ later acquires two square roots
of its own in $F$.  If $m\notin K$, then $z_5$ has no
square roots in $F$, and in no case does $z_5$ have any fourth roots in $F$.

Now if we build the automorphism
tree $I_F$ using the presentation $F=\{x_0=z_0,x_1=z_1,\ldots\}$ exactly as above,
then the level function $l_F$ is not computable:
$I_F$ contains the nodes $\la x_0, x_1\ra$ and $\la x_0,x_2\ra$,
and also contains $\la x_5\ra$, but it is impossible to decide
whether $\la x_5\ra$ has no immediate successors or two,
since this depends on whether $m\in K$.
(One would repeat this strategy, adding other field elements
to $F$ as roots of primes other than $2$, to code whether
each other $m'\in\omega$ lies in $K$ or not.)

On the other hand, now consider the distinct presentation
$E=\{ y_0=z_3,y_1=z_0, y_2=z_1,y_3=z_2,y_4=z_4,\ldots,y_{i}=z_i,\ldots\}$
of the same field $F$.  Using this presentation, we get a different sequence of
minimal polynomials:  $y_0$ has minimal polynomial $Y_0^8-2$,
so $q_0(Y_0)=Y_0^8-2$.  Then $y_1$ has minimal polynomial $Y_1-y_0^4$ over $\Q[y_0]$,
since $y_1=y_0^4$ is in this field.  That is, $q_1(Y_0,Y_1)=Y_1-Y_0^4$.
Next, $q_2(Y_0,Y_1,Y_2)= Y_2-Y_0^2$, since $y_2=z_1=z_3^2=y_0^2$,
and $q_3=Y_3+Y_2$, and $q_4=Y_4+Y_0$, and thereafter each $q_i$
is the same as in the previous presentation, except with the roles of $Y_0$ and $Y_3$
reversed.  This seemingly trivial difference between $E$ and $F$
changes the structure of the corresponding automorphism tree $I_E$:
now the only nodes at level $1$ are $\la y_0\ra$ and $\la y_4\ra$,
each with exactly one successor at level $2$,
exactly one at level $3$, at level $4$, at level $5$, and still
exactly one at level $6$ (these being $\la y_0,y_1, y_2,y_3,y_4,y_5\ra$
and $\la y_4,y_1, y_2,y_3,y_0,y_5\ra$), and so on.
So in this case the level function for $I_E$ is computable!
Since the eighth roots of $2$ appeared first, the question
of how many fourth roots of $2$ lie in $E$ is obviated.
It is possible, of course, that two more fourth roots of $2$
(specifically, square roots of $z_5$) may appear in $E$,
as $z_i$ and $z_j$ for some large $i<j$, but if so,
then we will simply have $q_i=X_i^2-X_5$ and $q_j=X_j+X_i$
for the presentation $F$, and likewise for the
presentation $E$, and we will be able to compute the number
of nodes at those levels in each presentation.  In $E$,
the appearance of $z_3$ before $z_0$ and $z_5$
allows the level function $l_E$ to evade the noncomputable
question of whether $z_5$ has square roots.

Of course, this micro-example is really the basic module
for the construction of two computable presentations $E$ and $F$
of a single field, with the property that $I_E$ has computable
level function but $I_F$ does not.  In the basic module above,
one would adjoin two square roots of $z_5$ when and if $0$ enters $K$.
Likewise, for every $m$, start with a full complement of $p_m$-th roots
of $2$, exactly one of which has $p_m$-th roots and also $(p_m^2)$-th
roots of its own.  If $m$ enters $K$, then adjoin all remaining $(p_m^2)$-th roots
of $2$ to the field.  The presentations differ exactly as in the basic module,
with the $(p_m^2)$-th roots of $2$ appearing before the $(p_m^3)$-th roots
of $2$ in the field $F$, but after them in the field $E$.
One can compute (for arbitrary $m$) the number of
$(p_m^3)$-th roots of $2$ in this field, but not the number of $(p_m^2)$-th roots.

\section{Orbit Relations on Fields}
\label{sec:Galoisaction}

In this section we wish to consider the action on a field $F$
of the automorphism group $\aut{F}$ of $F$.  We continue to assume that
$F$ is computable and algebraic.  Of course, the automorphism
group may have cardinality as high as $2^\omega$,
making it difficult to present in an effective fashion.
However, since the field $F$ may be viewed as
the union of an effectively presented chain of
finitely generated subfields, we will be able
to make substantial use of the following definition.

\begin{defn}
\label{defn:Galoisaction}
Let $F$ be any computable field, and let $G$
be any subgroup of the automorphism group of $F$.
The \emph{full action of $G$ on $F$} is the set
$$ \set{\la a_0,\ldots,a_{n-1},b_0,\ldots,b_{n-1}\ra}{
(\exists\sigma\in G)(\forall i < n)\sigma(a_i)=b_i},$$
where $\avec$ and $\bvec$ are tuples of elements of $F$,
with every $n$ allowed as their common length.
When $F$ is algebraic, we will be able to restrict
our attention to the \emph{action of $G$ on $F$},
which by definition is the set
$$ \set{\la a,b\ra}{a,b\in F~\&~(\exists\sigma\in G)\sigma(a)=b}.$$
If this set is computable, we say that $G$
\emph{acts computably on $F$}.

If $G=\aut{F}$, we also call these the \emph{full orbit relation}
and the \emph{orbit relation} of $F$, respectively.
\end{defn}
It is quite possible for an uncountable
$G$ to act computably on $F$:  for example,
the entire automorphism group $G$ of
(a computable presentation of) the algebraic closure
$\Qbar$ has size $2^\omega$, and contains
elements of every Turing degree, yet
its full action on $\Qbar$ is computable.
This will follow from Lemma \ref{lemma:fullaction}
and Corollary \ref{cor:normal} below.

Since we restrict our attention in this paper to algebraic
fields, we will only consider the action of a $G$ on an $F$,
not the full action.  The following lemma justifies this.
\begin{lemma}
\label{lemma:fullaction}
Let $F$ be a computable algebraic field of characteristic $0$,
and $G$ a subgroup of $\aut{F}$.
Set $B$ to be the action of $G$ on $F$,
and $A$ the full action.  Then $B\equiv_1 A$.
In particular, each is computable iff the other is,
and each is c.e.\ iff the other is.
\end{lemma}
In characteristic $p$, the proof below
ensures that $B\equiv_m A$.  $1$-equivalence
would be false in a finite field $F$.
\begin{pf}
$B\leq_1 A$ is immediate, so we build a function
$g$ to show $A\leq_1 B$.
Let $(a_0,\ldots,a_n;b_0,\ldots,b_n)\in F^{2n+2}$.
First we check whether the map $a_i\mapsto b_i$ extends
to an isomorphism from $Q(\avec)$ onto $Q(\bvec)$,
where $Q$ is the prime subfield of $F$.  This is computable,
since we need only find the minimal polynomial of each $a_i$
over $Q(a_0,\ldots,a_{i-1})$ and check that $b_i$ satisfies the
corresponding polynomial over $Q(b_0,\ldots,b_{i-1})$,
with each $a_j$ ($j<i$) mapped to $b_j$ to determine
the coefficients in this corresponding polynomial.
If $a_i\mapsto b_i$ does not extend to an isomorphism,
then clearly $\la\avec,\bvec\ra$ does not lie in the full action
of $G$ on $F$, and we define $g(\la\avec,\bvec\ra)$
to lie outside $B$.
(To ensure injectivity, let it be the $\la\avec,\bvec\ra$-th element
of some infinite c.e.\ subset of $\Bbar$.)
If it does extend to an isomorphism, we
use Theorem \ref{thm:primitive}
to find a single element $a\in Q(\avec)$ and polynomials
$p\in Q[X_0,\ldots,X_n]$ and $q_i\in Q[X]$ such that
$a=p(\avec)$ and $a_i=p_i(a)$ for each $i\leq n$.
(Since there are infinitely many such $a$, we may choose
ours to preserve injectivity of $g$.)
By the isomorphism above, each $b_i=p_i(p(\bvec))$ as well.
So $\la\avec,\bvec\ra$ lies in $A$ iff
$g(\la\avec,\bvec\ra)=\la a, p(\bvec)\ra$ lies in $B$.
\qed\end{pf}

On its face, membership in $B_F$ is a $\Sigma^1_1$ property:
it demands the existence of a function from $\omega$ to
$\omega$ satisfying certain arithmetic properties.  In fact, though,
the algebraicity of the field $F$ makes $B_F$ (and its
computable isomorph $A_F$) vastly simpler than this.
The proof demonstrates the usefulness of the isomorphism
trees defined in Section \ref{sec:isotrees}.
\begin{prop}
\label{prop:coce}
Let $F$ be a computable algebraic field.  Then
an arbitrary pair $\la a;b\ra$ of elements of $F$
lies in the field's orbit relation $B_F$ iff:
\begin{itemize}
\item
$b$ is conjugate to $a$ over the prime subfield $Q$; and
\item
$(\forall p\in Q[X,Y])[\text{$p(a,Y)$~has a root in $F$ $\implies$
$p(b,Y)$ has a root in $F$}]$.
\end{itemize}

It follows that the orbit relation $B_F$
of $F$ is $\Pi^0_2$.  More specifically, $B_F$ is $\Pi^S_1$,
where $S$ is the splitting set of $F$, and so if $F$ has a splitting
algorithm, then its orbit relation is co-c.e.
\end{prop}
By Lemma \ref{lemma:fullaction}, these results also apply to $A_F$, of course.
The surprise is that the second condition is not symmetric in $a$ and $b$;
this is essentially a consequence of Corollary \ref{cor:endo},
which can be used to show that the second condition is equivalent
to the same statement with $a$ and $b$ interchanged.

\begin{pf}
If $\la a,b\ra\in B_F$, then the two conditions given are immediate.
For the backwards direction,
let $\{y_0,y_1,\ldots\}$ be the domain of $F$.  (Of course,
this domain is really $\omega$; we write $y_n$ instead of $n$
to avoid confusion with the language of fields.)
Given $\la a,b\ra\in F\times F$, define $x_0=a$ and $x_{s+1}=y_s$
for all $n$.  So $\la x_s\ra_{s\in\omega}$ is a computable
generating sequence for $F$, and we may construct
the automorphism tree $I_F$ relative to this sequence,
along with the sequence $\la q_s\ra_{s\in\omega}$
of minimal polynomials, as in Definition \ref{defn:isotree}.
Thus $q_0(X)$ is the minimal polynomial of $a$
over the prime subfield $Q$.  Of course, some $n>0$ has
$x_n=a$, but this only means that $q_n(X_0,\ldots,X_n)=X_n-X_0$.

The first condition in the proposition is that $q_0(b)=0$.
Since we are assuming that both conditions hold,
$\la b\ra$ is a node at level $1$ in $I_F$,
and we claim that $\la b\ra$ lies on a path through $I_F$.
To see this, fix any $n\in\omega$ and consider a primitive generator $u\in F$
of the subfield $Q(a,x_1,x_2,\ldots,x_n)$.  Choose $p\in Q[X,Y]$
so that $p(a,Y)$ is the minimal polynomial of $u$ over $Q(a)$.
Since $F$ satisfies the second condition in the proposition,
$p(b,Y)$ must also have a root $v$ in $F$,
and so there is an isomorphism $h:Q(a,u)\to Q(b,v)$ with $h(a)=b$
and $h(u)=v$.  But then $\la b,h(x_1),h(x_2),\ldots, h(x_n)\ra\in I_F$,
by the definition of $I_F$.  Since this holds for all $n$, the node $\la b\ra$ has
arbitrarily long successors in $I_F$, and so,
applying K\"onig's Lemma to the finite-branching tree $I_F$,
$\la b\ra$ lies on a path through $I_F$.  This path defines an
automorphism of $F$ mapping $a$ to $b$, so $\la a,b\ra\in B_F$
as required.
\qed\end{pf}

Indeed, the preceding proof showed more than Proposition
\ref{prop:coce} stated.  We could have used any generating set for $F$
in place of the domain $\{ y_0,y_1,\ldots\}$, and constructed the automorphism
tree $I_F$ relative to this sequence (with $a$ attached as the first element
of the sequence).  The proof really showed that $\la a,b\ra\in B_F$
iff the node $\la b\ra$ had successors at arbitrary high levels in $I_F$.
\begin{cor}
\label{cor:coce}
Let $\{z_0,z_1,\ldots\}$ be a computable sequence of elements of
a computable algebraic field $F$ with prime subfield $Q$, such that
$$ Q\subseteq Q(z_0)\subseteq Q(z_1)\subseteq
\cdots~~~~\text{and}~~~~\cup_s Q(z_s)=F.$$
A pair $\la a,b\ra$ of elements of $F$ lies in the field's orbit relation $B_F$
iff $a$ and $b$ are conjugate over $Q$
and, for every $s$ and every $p_s\in Q[X,Y]$ such that $p_s(a,Y)$
is the minimal polynomial of $z_s$ over $Q(a)$, it holds that $p_s(b,Y)$
has a root in $F$.
\end{cor}
\begin{pf}
For each $s$, choose $y_s$ so that $p_s(b,y_s)=0$.
We apply Corollary \ref{cor:Konig} to $F$ and the subfield of $F$
generated by $\{b,y_0,y_1,\ldots\}$, with the common subfield
being $Q(a)$ (within $F$) and $Q(b)$ (within $Q(b,y_0,y_1,\ldots)$),
identified via the isomorphism mapping $a$ to $b$.
Thus these two fields are isomorphic over the common subfield,
so we have a field embedding of $F$ into itself mapping $a$ to $b$.
By Corollary \ref{cor:endo}, this embedding is an automorphism of $F$.
\qed\end{pf}

If $F$ is a normal algebraic extension of $Q$, then for any $a,b\in F$
which are conjugate over $Q$, both conditions in Proposition \ref{prop:coce} hold.
We state the obvious corollary.
\begin{cor}
\label{cor:normal}
All normal computable algebraic fields have computable orbit relation.
\qed\end{cor}

\section{Fields with Splitting Algorithms}
\label{sec:splitalg}

\begin{thm}
\label{thm:splitalg}
Let $F$ be a computable algebraic field with
a splitting algorithm.  Then 
$F$ is computably categorical iff
the orbit relation of $F$ (or equivalently,
the full orbit relation of $F$) is computable.

\end{thm}
\begin{pf}
We prove the forwards direction as Proposition \ref{prop:splitalg} below.
For the converse (which is the easier direction), we work under
the weaker assumption that $F$ has a computable presentation
with computable full orbit relation and with computable level function,
as in Definition \ref{defn:isotree}.  Of course, all computable algebraic
fields with splitting algorithms have computable level functions, so
this will suffice.  Since computable categoricity
is a property of the isomorphism type, we may assume $F$
itself to be the computable presentation which has a computable
level function $l_F$.
We take the domain of $F$ to be $\omega$.

Suppose the full orbit relation $A$ of $F$ is computable,
and let $\Ftilde\cong F$ be a computable copy of $F$.
We build a computable embedding $f:F\to\Ftilde$ as follows.
The prime subfield $F_0=Q$ is a computable subfield of $F$
(since $F$ is algebraic) and
has a unique embedding $f_0$ into $\Ftilde$,
which is computable and extends to the given isomorphism
$\theta_0$ (not necessarily computable) from $F$ onto $\Ftilde$.

We now proceed by recursion on $s$.
Given the embedding $f_s:F_s\embeds\Ftilde$
which we have already built,
consider the least element $z$ of $(F-F_s)$,
and let $F_{s+1}=F_s[z]$, on which we will define
$f_{s+1}$ to extend $f_s$.
By induction, $F_s$ must then be generated by
the finite set $\{ 0,1,\ldots,z-1\}$, and we also know inductively
that $f_s$ extends to some isomorphism $\theta_s$ from $F$ onto $\Ftilde$.
The splitting algorithm for $F_s$ lets us find the minimal
polynomial $q(Z)\in F_s[Z]$ of $z$.  Now the roots $z'\in F$
with $q(z')=0$ correspond precisely to the nodes
$\sigma\in I_F$ of the form $\la 0,1,\ldots,z-1,z'\ra$,
so let $d=l_F(\la 0,1,\ldots,z-1\ra)$ be the number of
such roots.  Then we can find all roots of $q(Z)$ in $F$:
let them be $z= z_{1} < z_{2}<\cdots < z_{d}$.
From the computable set $A$, we can determine
exactly which tuples $\xvec_i=\la 0,1,\ldots,(z-1),z, 0,1,\ldots, (z-1),z_{i}\ra$ lie in $A$;
these $z_i$ are the true conjugates of $z$ in $F$ over $F_s$.
For each $i\leq d$ with $\xvec_i\notin A$,
we know that the node $\sigma_i=\la 0,1,\ldots,z-1,z_i\ra$
is nonextendible in $I_F$, and so we use the function $l_F$
to find a level which contains no successors of $\sigma_i$.
(K\"onig's Lemma shows that such a level must exist.)
Thus we compute a single level $n$ such that no $\sigma_i$
with $\xvec_i\notin A$ has any successor at level $n$ in $I_F$.
Then we turn to the isomorphism tree $I_{F\Ftilde}$.
$\theta_s$ shows that the node $\tau=\la f_s(0),\ldots,f_s(z-1)\ra$
is extendible in $I_{F\Ftilde}$; indeed it is the image
of $\la 0,\ldots,z-1\ra$ under the isomorphism $H_{\theta_s}$
from Lemma \ref{lemma:autotree}.  Moreover,
from $H_{\theta_s}$ and our knowledge of $I_F$,
we know that any immediate successor of this $\tau$
which extends to level $n$ must be extendible.
So we enumerate $I_{F\Ftilde}$ until we find in it
some node $\rho$ at level $n$ extending $\tau$,
and we define $f_{s+1}(z)=\rho(z)$.  Now
$F_{s+1}=F_s[z]$, so $f_{s+1}$ is uniquely defined on $F_{s+1}$.
Thus $F_{s+1}$ is generated by $\{ 0,\ldots,z\}$, and moreover
$f_{s+1}$ extends to some isomorphism $\theta_{s+1}$
from $F$ onto $\Ftilde$, namely the (not necessarily unique)
path through $\rho\res (z+1)$ in $I_{F\Ftilde}$.
These were all the inductive facts we needed
in order to continue to the next stage.

This process computes $f_{s+1}$, uniformly in $s$, so
$f=\cup_s f_s$ is a computable embedding of $F$
into $\Ftilde$.  By Corollary \ref{cor:Konig}, it must be an isomorphism.
Thus $F$ is computably categorical.
\qed\end{pf}

To help the reader, we give a quick translation of the preceding
construction into the language of fields, without using
isomorphism trees. To find an image $f_{s+1}$ for $z$ in $\Ftilde$,
we find all conjugates of $z$ over $F_s$ in $F$, and use the
computable set $A_F$ to determine which are true conjugates.
For each false conjugate $z_i$, we search for a polynomial $p_i\in F_s[Z,Y]$
which shows $z_i$ to be false, namely, a polynomial such that
$p_i(z,Y)$ has a root in $F$ but $p_i(z_i,Y)$ does not.
Proposition \ref{prop:coce} shows that we will
eventually find such a polynomial.  (Of course, we are
using the computability of $S_F$ here to determine that
$p_i(z_i,Y)$ has no root in $F$.)
Using $f_s$, we then find the corresponding polynomials
$\ptilde_i\in (f_s(F_s))[Z,Y]$, and search for any $\ztilde\in\Ftilde$
which satisfies the (image in $(f_s(F_s))[X]$ of the) minimal polynomial
of $z$ over $F_s$ and such that every
$\ptilde_i(\ztilde,Y)$ has a root in $\Ftilde$.
We eventually must find such a $\ztilde$ (since
$\theta_s(z)$ is such an element), and when we do,
we know that it must be a true conjugate of $\theta_s(z)$,
since a false conjugate $\ztilde'$ would be the image of some false conjugate
$z_i$ of $z$ over $F_s$ in $F$, and therefore $\ptilde_i(\ztilde',Y)$
would not have had a root in $\Ftilde$.
So it is safe for us to define $f_{s+1}(z)=\ztilde$,
and $\theta_{s+1}$ is the composition of $\theta_s$
with an automorphism of $\Ftilde$ sending $\theta_s(z)$
to its true conjugate $\ztilde$ over the image $f_s(F_s)$.

\begin{prop}
\label{prop:splitalg}
If an algebraic computable field $F$ is
computably categorical, then the orbit relation $B_F$
of $F$ is computable.
\end{prop}
\begin{pf}
We will construct a computable field $\Ftilde$,
isomorphic to $F$, in such a way that the existence
of any computable isomorphism from $F$ onto $\Ftilde$
will allow us to compute $B_F$.  This is sufficient
to prove the proposition.  (It is also the contrapositive
of the usual argument for this sort of theorem.
In other contexts, mathematicians have often
taken the given property -- in this case,
the noncomputability of $B_F$ -- and used it to build
the second structure $\Ftilde$ by a construction which
diagonalizes against all possible computable
isomorphisms.  We believe that this would be feasible
in the present case, and that the construction would not
be substantially different from ours, but we see
our argument as more direct.)

We start by enumerating a generating set for the given field $F$.
Set $F_0$ to be the prime subfield of $F$, either $\Q$ or $\F_p$,
with $z_0$ as the multiplicative identity element of $F$.
Given $F_s$, choose the least number $y\in F-F_s$, and
let $F_{s+1}$ be the normal closure of $F_s\cup\{ y\}$ within $F$.
Since $F$ has a splitting algorithm, this is computable:
we can determine the minimal polynomial $p(X)\in F_0[X]$
of $y$ over $F_0$ and then find all of its roots in $F$,
and $F_{s+1}$ is generated by these roots over $F_s$.
Being finitely generated over $F_s$ (and hence over $F_0$,
by induction), $F_{s+1}$ has a splitting algorithm,
uniformly in its generators and hence uniformly in $s$,
and thus is a decidable subset of the algebraic field $F$.
We set $z_{s+1}$ to be the least primitive generator
of $F_{s+1}$ over $F_0$, and list out all of its
conjugates over $F_0$ in $F$, each of which is another
primitive generator of $F_{s+1}$:
$$ z_{s+1}=z_{s+1}^0 < z_{s+1}^1 <\cdots < z_{s+1}^{d_{s+1}}.$$
(The superscripts here are not exponents, of course,
but merely indices.)
Again, all of this is computable uniformly in $s$.
Notice that $$d_{s+1}+1 =[F_{s+1}:F_0]$$ is precisely
the size of the automorphism group of $F_{s+1}$;
this is at most the degree of the minimal polynomial $p_{s+1}(X)\in F_0[X]$
of $z_{s+1}$ over $F_0$, and in general is not equal
to that degree, since $F_{s+1}$ may fail to be normal over $F_0$.
However, every automorphism $\sigma$ of $F_{s+1}$
must fix setwise (although not necessarily pointwise)
every $F_t$ with $t\leq s$, since each $F_t$ is normal within $F$,
hence normal within $F_{s+1}$.  (In particular, $\sigma(z_t)$
must equal some $F_0$-conjugate of $z_t$, which is to say,
some primitive generator $z_t^i$ of $F_t$.)

The requirements for our construction are simply stated:
$$\mR_e :~~\phi_e\text{~is not an isomorphism from $F$
onto~}\Ftilde.$$
If all these requirements (for all $e$) were true of the field
$\Ftilde$, then there would be no computable isomorphism from
$F$ onto $\Ftilde$.  Our construction
of $\Ftilde$ addresses these requirements individually, for each $e$,
and attempts to satisfy each one.
Requirements such as these are used throughout computability
theory; the reader unfamiliar with them should consult \cite{S87}.
Normally, to prove computable non-categoricity of $F$, one
would build a computable field $\Ftilde$ isomorphic to $F$
for which every $\mR_e$ holds.  As we are proving the contrapositive,
our $F$ is assumed to be computably categorical, and so we will not succeed
in satisfying all of these requirements, but our construction will
attempt to do so nevertheless, using the indices $e$
to assign priorities to each requirement, with a lower index
denoting a higher priority.  The least $e$ for which we fail to
satisfy $\mR_e$ will be the key to our decision procedure for $B_F$.

Now we construct $\Ftilde$, in stages, with each $\Ftilde_s$
isomorphic to $F_s$.  Of course $\Ftilde_0=F_0$ is just a computable
copy of the prime subfield, and we fix $$f_0:F_0\to\Ftilde_0$$
to be the (unique) isomorphism between them.
All requirements $\mR_e$ are unsatisfied at stage $0$.
At the end of each stage $s$, we will have a
field $\Ftilde_s$ and an isomorphism $f_s$ onto it from $F_s$,
all computable uniformly in $s$.  Moreover, for each $t\leq s$,
we also know $z_t$ and all $z_t^i$, for $i\leq d_t$.
We proceed as follows to build $\Ftilde_{s+1}$.

For each $t\leq s$ and $i\leq d_t$, we consider the automorphism
$\sigma_t^i$ of $F_t$ defined by $\sigma_t^i(z_t)=z_t^i$;
these are precisely the automorphisms of $F_t$ (over $F_0$,
which is rigid).  We do the same for the automorphisms
$\sigma_{s+1}^k$ of $F_{s+1}$, defined by $\sigma_{s+1}^k(z_{s+1})=z_{s+1}^k$,
for $k\leq d_{s+1}$.  Since all this is computable, we
may decide, for each $t$ and $i$, whether $\sigma_i^t$ extends
to an automorphism of $F_s$ and/or to an automorphism of $F_{s+1}$.
Certainly $\sigma_t^0$, the identity, extends to
$\sigma_{s+1}^0$, but other $\sigma_t^i$ may or may not extend
to $F_s$, and those which do may or may not extend to $F_{s+1}$.
However, by normality of $F_s$ within $F_{s+1}$, $\sigma_t^i$ can only
extend to an automorphism of $F_{s+1}$ if it extends to an automorphism of $F_s$.

We search for the least $e\leq s$, if any, for which
$\mR_e$ is not yet satisfied and there exists some $t\leq s$ and $j \leq d_t$
for which $\phi_{e,s}(z_t^i)\converges$
for all $i\leq d_t$
and $\sigma_t^j$ extends to an automorphism $\sigma$ of $F_s$
but does not extend to any automorphism of $F_{s+1}$.
(Our conventions about $\phi_{e,s}(z_t^i)$ and $W_{e,s}$
were described in Section \ref{sec:computability}.)
If there is no such $e$, then we define
$f_{s+1}\res F_s =f_s$ and let $\Ftilde_{s+1}$
contain all of $\Ftilde_s$, along with fresh elements
to be the images of the elements
of $(F_{s+1}-F_s)$ under $f_{s+1}$.

If such an $e$ does exist, then we act to satisfy $\mR_e$,
using the value $j$ and the automorphism $\sigma$ of $F_s$ found above.
First, though, for all of the (finitely many) elements $x\in W_{e,s}\cap F_t$,
we find a polynomial $q(Z)\in F_0[Z]$ with $q(z_t)=x$
and check whether $q(\phi_e(z_t))=\phi_e(x)$.
If this fails for any $x$, then $\phi_e$ cannot be an isomorphism,
so we act just as we did (above)
when $e$ did not exist, and declare $\mR_e$ satisfied.
Otherwise, we know that $\phi_{e,s}$ maps the set
$\{ z_t^0,\ldots,z_t^{d_t}\}$ bijectively onto the set
$\{ f_s(z_t^0),\ldots, f_s(z_t^{d_t})\}$.
Fix the $m$ such that $f_s(z_t^m)=\phi_e(z_t^0)$.
If no automorphism of $F_s$ maps $z_t^0$ to $z_t^m$,
then $(f_s^{-1}\circ\phi_e)$ cannot be an automorphism,
so $\phi_e$ cannot be an isomorphism,
and once again we just extend $f_s$ to $f_{s+1}$, fill in
$\Ftilde_{s+1}$ with fresh elements, and declare $\mR_e$ satisfied.
Otherwise, there is an automorphism $\tau$ of $F_s$ with $\tau(z_t^0)=z_t^m$,
and we let $f_{s+1}\res F_s = f_s\circ \tau\circ \sigma^{-1}$,
and form $\Ftilde_{s+1}$ by adding fresh elements to $\Ftilde_s$
to be the images of the elements of $F_{s+1}$ under $f_{s+1}$.
Lemma \ref{lemma:satisfy} below will show that in this case,
$\phi_e$ cannot be an isomorphism from $F$ onto $\Ftilde$.
So we declare $\mR_e$ satisfied, and end the stage.

This builds a computable field $\Ftilde$, which we claim is isomorphic
to $F$.  Of course, we made no attempt during the construction to ensure
that $\lim_s f_s(x)$ must exist for $x\in F$.  However,
every finitely generated subfield of $F$ embeds into some $F_s$ and hence
(via $f_s$) into $\Ftilde_s$, whence into $\Ftilde$.  A symmetric argument
with $\Ftilde$ and $F$ interchanged also holds, leaving Corollary \ref{cor:Konig}
to prove that $F\cong\Ftilde$ over $F_0$.  Of course, the isomorphism
need not be computable.

\begin{lemma}
\label{lemma:satisfy}
If $\mR_e$ is ever declared satisfied during this construction,
then $\phi_e$ is not an isomorphism from $F$ onto $\Ftilde$.
\end{lemma}
\begin{pf}
Let $s+1$ be the stage at which $\mR_e$ is declared satisfied.
If $\phi_e$ were an isomorphism
from $F$ onto $\Ftilde$, then it would have
to restrict to an isomorphism from $F_{s+1}$ onto $\Ftilde_{s+1}$,
because $F_{s+1}$ is normal within $F$ and
$\Ftilde_{s+1}$ is the image of $F_{s+1}$ under the isomorphism $f$.
Similarly, $\phi_e\res F_s$ would have to be an isomorphism.
To prove the lemma, therefore, we will show that $\phi_e\res F_s$
and $\phi_e\res F_{s+1}$ cannot both be isomorphisms.

First, if any $\phi_{e,s}(x)$ converged to a value other than a root
of $q(\phi_e(z_t))$ (where $q(Z)$ is the minimal polynomial of $x$ over $F_0$,
as in the construction), then clearly
$\phi_e\res F_s$ is not an isomorphism.  Also,
we chose $m$ to satisfy $\phi_e(z_t^0)=f_s(z_t^m)$,
so if no automorphism of $F_s$
maps $z_t^0$ to $z_t^m$, then $(f_s^{-1}\circ\phi_e)\res F_s$
cannot be an automorphism, and $\phi_e\res F_e$ cannot be an isomorphism

In the remaining case, 
$\sigma_s^j$ is known to extend to an automorphism
$\sigma$ of $F_s$, but not to any automorphism of $F_{s+1}$,
and we have an automorphism $\tau$ of $F_s$ with $\tau(z_t^0)=z_t^m$.
In this case
$$ f_{s+1}(z_t^j) = (f_s\circ \tau\circ (\sigma^{-1}))(z_t^j)
= f_s(z_t^m) = \phi_e(z_t^0).$$
With the new elements added to form $\Ftilde_{s+1}$,
this means that $\phi_e$ cannot restrict to an isomorphism from
$F_{s+1}$ onto $\Ftilde_{s+1}$, because if it were,
then $(f_{s+1}^{-1}\circ\phi_e)$ would be an automorphism
of $F_{s+1}$ mapping $z_t^0$ to $z_t^j$, which would necessarily
extend $\sigma_s^j$, and we chose $j$ precisely so that
no such extension of $\sigma_s^j$ exists.
\qed\end{pf}

For the sake of readers who saw ``$\mR_e$'' and expected an injury
construction, we note that no conflict exists between the different
requirements here:  Lemma \ref{lemma:satisfy} shows that there was no need
to preserve the satisfaction of $\mR_e$ once it was established.
The normality of each $F_{s+1}$ within $F$ took care of that.
Moreover, the above argument did not require that $\lim_s f_s$
itself exist, let alone that it be an isomorphism, and so, when
satisfying a requirement, the construction makes no effort
to have $f_{s+1}$ agree with $f_s$ at all.
Therefore, there are no injuries in this construction.
The only need for priority arose in choosing which $\mR_e$ to satisfy,
at a stage at which more than one requirement might have been
satisfiable, and the only reason for taking the least $e$ in those situations
was that this is the simplest way to ensure that each
requirement which can be satisfied at infinitely many stages
does indeed eventually become satisfied.

Of course, it remains to show that this field $\Ftilde$ really does prove
the desired result.  We have seen that $\Ftilde\cong F$, and by the computable
categoricity of $F$, this means that there exists a computable isomorphism
$\phi_e$ from $F$ onto $\Ftilde$.  We fix the least such index $e$
and the largest stage $s_0$ at which a requirement $\mR_i$ with $i<e$
was satisfied.  By Lemma \ref{lemma:satisfy}, $\mR_e$ was never satisfied,
so at all stages $s>s_0$, it did not fulfill the conditions which the construction
posed in order to be satisfied.

We now describe an algorithm for deciding whether a pair
$\la z_t,z_t^n\ra$, with $t$ and $n\leq d_t$ arbitrary, lies in the orbit relation
$B_F$ of $F$.  We may assume that $t > s_0$, since the answers
for all $t\leq s_0$ and $n\leq d_t$ constitute finitely much information.
Find the least $s_1$ such that $\phi_{e,s_1}(z_t^i)\converges$ for all $i\leq d_t$.
Now $z_t^0\geq t-1$ (since $z_t^0\notin F_{t-1}$
and, by construction, each subfield $F_s$
contains the elements $0,\ldots,s-1$), and hence $\phi_e(z_t)$ requires
at least $t-1$ steps to converge, forcing $s_1 \geq t-1\geq s_0$.
As above, we may compute all automorphisms of $F_{s_1}$ and check
whether any of them maps $z_t$ to $z_t^n$.  If not, then clearly
$\la z_t,z_t^n\ra\notin B_F$, since any automorphism mapping $z_t$
to $z_t^n$ would restrict to an automorphism of the
subfield $F_{s_1}$, it being normal within $F$.  We claim
that if some automorphism $\rho$ of $F_{s_1}$ does have
$\rho(z_t)=z_t^n$, then $\la z_t,z_t^n\ra\in B_F$.

To see this, we induct on stages $s\geq s_1$,
claiming that $\sigma_t^n$ extends to an automorphism
of $F_s$ for all such $s$.  For $s_1$ this already holds,
since $$\rho(z_t)=z_t^n=\sigma_t^n(z_t)$$ and
$z_t$ generates the domain $F_s$ of $\sigma_t^n$.
So we consider an arbitrary $s+1 > s_1$.
Now since $\mR_e$ does \emph{not} become satisfied
at stage $s+1$ (and no $e' <e$ ever becomes satisfied after stage $s_0$),
$e$ must not fulfill the conditions in the construction for choosing
the requirement to be satisfied at stage $s+1$.
But with $s\geq s_1$, $\phi_{e,s}(z_t^i)\converges$
for all $i\leq d_t$, and $\sigma_t^n$ extends
(by inductive hypothesis) to an automorphism
of $F_s$.  If $\sigma_t^n$ failed to extend to an automorphism
of $F_{s+1}$, then the construction would have chosen $e$
and acted to satisfy $\mR_e$ at this stage, destroying the
isomorphism $\phi_e$.  This did not happen,
so $\sigma_t^n$ must extend to an automorphism of $F_{s+1}$.
This completes the induction.

The extension of $\sigma_t^n$ to an automorphism
$\tau$ of the whole field $F$
is now accomplished by application to Corollary
\ref{cor:Konig}.  We let $Q=F_t$, as a subfield of $F$,
and have $E=F$.  The subfield $Q$ of $E$ is in fact
$F_t$ as well, but we identify the two copies
of $F_t$ via $\sigma_t^n$, rather than via the identity.
Corollary \ref{cor:Konig} then yields an isomorphism
between $F$ and $E$ over the common subfield,
i.e.\ an automorphism of $F$ extending $\sigma_t^n$,
as desired.

Finally, for an arbitrary pair $\la a,b\ra$ of elements of $F$,
we find an $s$ with $a,b\in F_s$ and determine all automorphisms
(if any) of $F_s$ mapping $a$ to $b$.  Each of these
automorphisms is equal to $\sigma_s^i$ for some $i\leq d_s$,
and so for each such $i$, we check whether $\la z_s,z_s^i\ra\in B_F$.
If this holds for any $i\leq d_s$, then also $\la a,b\ra\in B_F$, as witnessed
by the automorphism(s) of $F$ extending $\sigma_s^i$.
Conversely, any automorphism of $F$ mapping $a$ to $b$
would restrict to an automorphism of $F_s$, which would then equal
$\sigma_s^i$ for one of these $i$.  Thus $B_F$ is computable.
\qed\end{pf}



\section{Relativizing the Results}
\label{sec:relative}

For simplicity, we proved Theorem \ref{thm:splitalg}
above in a non-relativized form.  However, the argument
in one direction relativizes easily to any degree $\bfd$,
producing the following result.  Recall that a
computable structure $\A$ is \emph{$\bfd$-computably categorical}
if every computable structure $\B$ classically
isomorphic to $\A$ is $\bfd$-computably isomoprhic to $\A$.

\begin{prop}
\label{prop:dcc}
Let $F$ be a computable algebraic field with
a splitting algorithm, and fix any Turing degree $\bfd$.
If the orbit relation of $F$ (or equivalently,
the full orbit relation of $F$) is $\bfd$-computable,
then $F$ is $\bfd$-computably categorical.
\qed\end{prop}

The proof of Proposition \ref{prop:splitalg},
however, relativizes to the statement that,
for each computable algebraic field $F$
with a splitting algorithm,
there exists a $\bfd$-computable field $\Ftilde\cong F$
such that any $\bfd$-computable isomorphism
from $F$ onto $\Ftilde$ would allow one to
compute $B_F$.  This is not sufficient to prove
the converse of Proposition \ref{prop:dcc}, and in fact
the converse turns out to be false, by the following argument.

By \cite{M09b}, there exists a computable
algebraic field $F$, with a splitting algorithm, such that
there is no least Turing degree $\bfd$ for which
$F$ is $\bfd$-computably categorical.
Indeed, it is proven there that there
exist Turing degrees $\bfc$ and $\bfd$
whose infimum is $\bfz$, such that $F$ is both
$\bfc$-computably categorical and $\bfd$-computably categorical.
If the orbit relation $B_F$ for this field
were both $\bfc$-computable and $\bfd$-computable,
then it would be computable, and $F$ would have been
computably categorical.  Consequently,
one of these degrees (say $\bfd$) has the property
that $\deg{B_F}\not\leq_T\bfd$, even though
$F$ is $\bfd$-computably categorical.
So the converse of Proposition \ref{prop:dcc} fails.
It would be of interest to determine
whether perhaps there exists a computable field $E$
isomorphic to this $F$, for which $\bfd$
does compute $B_E$ (in which case necessarily
$\bfc$ would not compute $B_E$, by the argument above).
If so, then the Turing degree of the orbit relation
would not be invariant under isomorphisms between
computable algebraic fields with splitting algorithms,
even though (by Theorem \ref{thm:splitalg})
the decidability of the orbit relation is invariant.

\section{Further Notes}
\label{sec:notes}

The article \cite{HKMS11} is complementary to this one in
a number of respects.  It considers relative computable categoricity
for computable algebraic fields, and also examines the possible
computable dimensions of such fields.  Its most relevant
results for us, however, are negative ones:
it is shown in \cite[Theorem 4.5]{HKMS11} that there exists a computable
algebraic field with computable orbit relation which is not computably
categorical, and it is shown in \cite[Theorem 5.1]{HKMS11} that there exists
a computably categorical algebraic field $F$ such that
$B_F$ is not even $\Sigma^0_2$, let alone computable.
One might have hoped for Proposition \ref{prop:splitalg}
to generalize to all computable algebraic fields;
alternatively, one might have rephrased Theorem
\ref{thm:splitalg} to say that computable categoricity
is equivalent to computable enumerability of $B_F$
(which is exactly the content of the proof, $B_F$ being
$\Pi^0_1$ for any computable field with a splitting algorithm).
The results in \cite{HKMS11} dash these hopes,
establishing that both directions of the desired equivalence
are false.

It is noted in \cite{HKMS11} that the definition of computable
categoricity, which is normally of complexity $\Pi^1_1$,
drops to complexity $\Pi^0_4$ when one restricts
the discussion to computable algebraic fields.
Essentially this follows from Corollary \ref{cor:Konig}
above (with $Q$ as the prime subfield), which reduces
the complexity of the isomorphism relation on such fields
dramatically.  Moreover, \cite[Theorem 6.4]{HKMS11}
proves that for algebraic fields, computable categoricity
is $\Pi^0_4$-complete, and thus is quantifiably more difficult
than computable categoricity for algebraic fields with
splitting algorithms.  This re-establishes the negative
results from the preceding paragraph.  Also, the proof
in \cite{HKMS11} that computable categoricity does not imply
relative computable categoricity is already of interest,
since to our knowledge, all previous results giving
structural criteria for computable categoricity in
commonplace mathematical classes also implied
relative computable categoricity.  So there is concrete evidence
that fields, and even just algebraic fields, constitute
a more challenging class of structures for this question
than did the previous classes studied.

\parbox{4.7in}{
{\sc
\noindent
Department of Mathematics \hfill \\
\hspace*{.1in}  Queens College -- C.U.N.Y. \hfill \\
\hspace*{.2in}  65-30 Kissena Blvd. \hfill \\
\hspace*{.3in}  Flushing, New York  11367 U.S.A. \hfill \\
Ph.D. Programs in Mathematics \& Computer Science \hfill \\
\hspace*{.1in}  C.U.N.Y.\ Graduate Center\hfill \\
\hspace*{.2in}  365 Fifth Avenue \hfill \\
\hspace*{.3in}  New York, New York  10016 U.S.A. \hfill}\\
\medskip
\hspace*{.045in} {\it E-mail: }
\texttt{Russell.Miller\at {qc.cuny.edu} }\hfill \\
}

\parbox{4.7in}{
{\sc
\noindent
East Carolina University\hfill \\
\hspace*{.1in}  Department of Mathematics \hfill \\
\hspace*{.2in}  Greenville, NC 27858 U.S.A. \hfill}\\
\medskip
\hspace*{.045in} {\it E-mail: }
\texttt{shlapentokha\at {ecu.edu} }\hfill \\
}

\end{document}